\documentclass[12pt]{article}
\usepackage{hyperref}

\newcommand{\RNum}[1]{\uppercase\expandafter{\romannumeral #1\relax}}
\usepackage{amsmath}
\usepackage{hyperref}
\usepackage{amssymb}
\usepackage{latexsym}
\usepackage{epsf}
\usepackage{supertabular}
\usepackage{rotating}
\usepackage{multicol}
\usepackage{multirow}
\usepackage{epsfig}
\usepackage{fancyhdr}       % fuer Deckblatt
\usepackage{booktabs}
\usepackage{graphics}
\usepackage{array}
\usepackage{cite}
\usepackage[numbers,sort&compress]{natbib}
\usepackage{cleveref}
\usepackage{makecell}

\newif\ifger
\gerfalse

\newtheorem{theorem}{Theorem}[section]

\newtheorem{lemma}[theorem]{Lemma}

\newtheorem{corollary}[theorem]{Corollary}

\newtheorem{proposition}[theorem]{Proposition}

\title{On symmetric 2-designs of prime order with almost simple flag-transitive automorphism groups
}

\author{ Ziwei Lu, Shenglin Zhou\footnote{Corresponding author. This work is supported by the National Natural Science Foundation of China (Grant No.12271173). slzhou@scut.edu.cn}\\
	\small \it School of Mathematics, South China University of Technology,\\
	\small \it Guangzhou 510641,  China}

\begin{document}
	\baselineskip=19pt
	\maketitle
	\date

\begin{abstract}
\indent In this article, we investigate symmetric $2$-designs of prime order admitting a flag-transitive automorphism group $G$. Recently, the authors proved that the automorphism group $G$ of this type of designs must be point-primitive, and is of affine or almost simple type. Here, we give the complete classification of symmetric $2$-designs of prime order, admitting a flag-transitive almost simple automorphism group.

		
\end{abstract}
	
\section{Introduction}
A \emph{symmetric $2$-$(v,k,\lambda)$ design} $\mathcal{D}=(\mathcal{P},\mathcal{B})$ is an incidence structure of $v$ points and $v$ blocks such that every point is incident with exactly $k$ blocks, and every pair of points is incident with $\lambda$ blocks. Points and blocks are interchangeable in this case, due to their dual role. The integers $v,k,\lambda$ called as the parameters of $\mathcal{D}$. We say that the design $\mathcal{D}$ is \emph{non-trivial} if $2<k<v-1$. If $\mathcal{D}=(\mathcal{P},\mathcal{B})$ is a symmetric 2-$(v,k,\lambda)$ design, then $\mathcal{D}^c=(\mathcal{P},\mathcal{B}^c)$ where $\mathcal{B}^c=\{\mathcal{P}\setminus B|B\in \mathcal{B}\}$ is a symmetric $(v,v-k,v-2k+\lambda)$ design which is called the \emph{complement} of $\mathcal{D}$. The value $n=k-\lambda$ is called \emph{the order of $\mathcal{D}$}, and clearly, $\mathcal{D}$ and $\mathcal{D}^c$ have the same order. An automorphism of a design $\mathcal{D}$ is a permutation of the points which also permutes the blocks, preserving the incidence relation. The set of automorphisms of a design with the composition of functions is a group, denoted by $Aut(\mathcal{D}).$ A \emph{flag} of $\mathcal{D}$ is an incident point-block pair $(\alpha,B)$. For a subgroup $G$ of $Aut(G)$, $G$ is said to be \emph{flag-transitive} (or \emph{point-primitive}) if $G$ acts transitively (or primitively) on the flags (or points) of $\mathcal{D}$. Note that the complement of a point-primitive design is
also point-primitive with the same full automorphism group.\\
\indent	A group $G$ is said to be \emph{almost simple} with socle $X$ if $X\unlhd G\le Aut(X),$ where $X$ is a non-abelian simple group. In \cite{R1}, Regueiro firstly reduced a flag-transitive point-primitive biplane must be of affine type or almost simple type. Later, she thoroughly completed the classification of biplanes with flag-transitive automorphism groups of almost simple type in \cite{R2,R3,R4}. Some results on symmetric 2-designs suggests studying this type of designs admitting the flag-transitive automorphism group whose the socle is a non-abelian simple group, see for example \cite{lprime,triplane}. Recently, we proved that if $\mathcal{D}=(\mathcal{P},\mathcal{B})$ is a flag-transitive symmetric 2-design of prime order, then the automorphism group $G$ is point-primitive of affine type or point-primitive of almost simple type. This paper is a continuation of \cite{lu}. The aim of this article is to classify symmetric $2$-designs of prime order, admitting a flag-transitive almost simple automorphism group. We now state the main result as follows:
\begin{theorem}\label{th1}
Let $\mathcal{D}=(\mathcal{P},\mathcal{B})$ be a symmetric $2$-design of prime order $n$, admitting a flag-transitive automorphism group $G$ such that $X=Soc(G)$ is a non-abelian finite simple group. Then
one of the following holds:
\begin{enumerate}
	\item[\rm (1)] $\mathcal{D}$ is a $2$-$(11,5,2)$ design or its complement, and $G=PSL_{2}(11);$
	\item[\rm (2)] $\mathcal{D}$ is a projective plane of order $n$ or its complement, and $PSL_{3}(n)\unlhd G\le P\Gamma L_{3}(n)$. 
\end{enumerate}
\end{theorem}
By the proof of Theorem \ref{th1} and the classification of $2$-transitive symmetric designs \cite{kantor}, we have the following corollary.
\begin{corollary}\label{cor1}
	If $\mathcal{D}$ is a symmetric $2$-design of prime order, admitting an almost simple automorphism group $G$, then $G$ is flag-transitive if and only if $2$-transitive.
\end{corollary}

\indent In Section $2$ we list some numerical conditions for flag-transitive $2$-designs. In Section $3$ we show that Theorem \ref{th1}. In order to prove Theorem \ref{th1}, we observe that the result for the case where the socle is sporadic simple groups or exceptional groups of Lie type follows from \cite{sporadic,exceptional}. Hence, we mainly consider that the socle $X$ is an alternating group or a finite classsical simple group.
\section{Preliminary results}
We note that the condition  $n=k-\lambda$ is prime implies that $(k,\lambda)=1$ or $n$. If $(k,\lambda)=1$, then from \cite{1} we conclude that either $\mathcal{D}$ is a $2$-$(11,5,2)$ design and $G=PSL_{2}(11)$ or $\mathcal{D}$ is a projective plane $PG_{2}(n)$  where $n$ is prime. Therefore, we will only focus on $(k,\lambda)=n$ in the following. Let $k=nk^*$ and $\lambda=n\lambda^*$, where $k^*=\frac{k}{(k,\lambda)}, \lambda^*=\frac{\lambda}{(k,\lambda)}$. We firstly state some basic results on designs.
\begin{lemma}\label{l1}
	Let $\mathcal{D}$ be a symmetric $2$-$(v, k, \lambda)$ design and $G$ be a flag-transitive automorphism group of $\mathcal{D}$. Then the following hold:
	\begin{enumerate}
		\item[\rm (1)] $(v-k)\lambda=(k-1)(k-\lambda)$;
		\item[\rm (2)] $\lambda(v-1)=k(k-1)$. In particular, if the order of $\mathcal{D}$ is $n$ and $(k,\lambda)=n,$ then
		\begin{enumerate}
		\item[\rm (i)] $v=\dfrac{k(k-1)}{\lambda}+1=n(k^*+1)+\dfrac{n-1}{k^*-1}$, where $k=nk^*;$
		\item[\rm (ii)] $v\le2k-1<2k$;
		\end{enumerate}		
		\item [\rm (3)] 
		$k \mid |M|$,  where $M:=G_{\alpha}$ is the point stabilizer of $\alpha$ in $G$;\item [\rm (4)]
		If $d$ is any non-trivial subdegree of $G$, then $k \mid \lambda d$. In addition, $k^* \mid d$ and $k^* \mid (v-1,d)$. 
	\end{enumerate}
\end{lemma}
{\bfseries Proof}. Parts (1) and (2)(i) are see for example \cite[Proposition 1.1]{Osifodunrin}. Part (4) follows from \cite{Da}. (2)(ii) Note that $(k,\lambda)=n,$ then $(k-\lambda,\lambda)=(k,\lambda)=k-\lambda$. It follows that $k-\lambda\mid \lambda$ and $k-\lambda\le \lambda$. Combining it with Part (1), we get $v-k\le k-1$, that is, $v\le2k-1<2k$. (3) Since $G$ is flag-tranistive, $M$ is transitive on the blocks containing a point $\alpha$, Thus we have $k \mid |M|$. \vspace{3mm}$\hfill\square$\\
\indent The following famous Bruck-Ryser-Chowla theorem gives a necessary condition for the existence of symmetric designs.
\begin{lemma}\label{l2}  {\rm (Bruck-Ryser-Chowla Theorem \cite{design})}
	Let $v,k,\lambda$ be integers with $\lambda(v-1)=k(k-1)$ for which there exists a symmetric $2$-$(v,k,\lambda)$ design.
	\begin{enumerate}
		\item [\rm (1)]If $v$ is even, then $n=k-\lambda$ is a square.
		\item[\rm (2)]If $v$ is odd, then the equation $(k-\lambda)x^2+(-1)^\frac{v-1}{2}\lambda y^2=z^2$ has a solution in integers $x,y,z$ not all zero.
	\end{enumerate}
\end{lemma}

\indent From Lemma \ref{l2}, we know if the order $n=k-\lambda$ is prime, then $v$ is odd. As a consequence, we can simplify the proof of main theorem in \cite{lu} by just deal with the product type. The next lemma is an elementary result on almost simple groups.
\begin{lemma}\label{l3}
Let $G$ be an almost simple group with the socle $X$, and let $M$ be maximal in $G$ not containing $X$. Then $G=MX$ and $|M|\mid |Out(X)|\cdot |M\cap X|$.
\end{lemma}
{\bfseries Proof}. The proof is straightforward. Since $M$ is maximal in $G$ not containing $X$, then $G=MX$. Note that $G$ is an almost simple group with the socle $X$, that is, $X\unlhd G\le Aut(X),$ where $X$ is a non-abelian simple group. It follows that  $|G|$ divides $|Out(X)|\cdot |X|.$ Therefore, $|M|\mid |Out(X)|\cdot |M\cap X|.$\vspace{3mm}$\hfill\square$\\
\indent In this article, we assume that $q=p^f$, where $p$ is a prime number and $f$ is a positive integer. Given a positive integer $a$, let $p$ be a prime factor of $a$, we denote the $p$-part of $a$ by $a_{p}.$ In other words, $a_{p}=p^t$ with $p^t\mid a$ but $p^{t+1}\nmid a.$ In addition, $a_{p'}$ denotes the $p'$-part of $a$.
\begin{lemma}\label{l4}{\rm (\cite {irr})}
	If $X$ is a group of Lie type in characteristic $p$, acting on the set of cosets of a maximal parabolic subgroup, and $X$ is neither $PSL_{m}(q)$, $P\Omega_{2m}^{+}(q)$ $($with $m$ odd$)$, nor $E_{6}(q)$, then there is a unique subdegree which is a power of $p$.
\end{lemma}

\indent We remark that even in the cases excluded in Lemma \ref{l4}, many of the maximal parabolic subgroups still have the property as asserted, see \cite[Lemma 2.6]{linear}. In particular, for an almost simple group $G$ with $X=Soc(G)=E_{6}(q)$, if $G$ contains a graph automorphism or $M=P_{i}$ with $i$ one of $2$ and $4$, the conclusion of Lemma \ref{l4} is still true.\\
\indent In \cite{odd}, the classification of odd degree primitive permutation groups have been completed. Since $v$ is odd in this paper, then the next two lemmas play an important role in proof of Theorem \ref{th1}.
\begin{lemma}\label{l5}{\rm (\cite[Theorem] {odd})}
	Let $G$ be a primitive permutation group of odd degree $v$ on a set $\Omega$. Assume that the socle $X = X(q)$ of $G$ is a classical simple group with a natural projective module $V=V(m,q)$ where $q=p^f$ and $p$ prime, and let $M:=G_{\alpha}$ be the stabilizer of $\alpha\in \Omega$, then one of the following holds: 
	\begin{enumerate}
		\item [\rm (a)]if $q$ is even, then $M\cap X$ is a parabolic subgroup of $X$.
		\item [\rm (b)]if $q$ is odd then one of $\rm(i)$, $\rm(ii)$ below holds:	
		\begin{enumerate}
		\item [\rm (i)]$M=N_{G}(X(q_{0}))$ where $q=q_{0}^t$ and $t$ is an odd prime; 
		\item [\rm (ii)]$X$ is a classical group with natural projective module $V=V(m,q)$ and one of $(1)$-$(7)$ below holds: 
	      \begin{enumerate}
		\item [\rm (1)]$M$ is the stabilizer of a non-singular subspace $($any subspace for $PSL_{m}(q)$$)$;
		\item [\rm (2)]$M\cap X$ is the stabilizer of an orthogonal decomposition $V=\oplus V_{j}$ with all $V_{j}$'s isometric $($any decomposition $V=\oplus V_{j}$ with dim$(V_{j})$ constant for $X=PSL_{m}(q)$$)$;
		\item [\rm (3)]$X=PSL_{m}(q)$, $M$ is the stabilizer of a pair $\{U,W\}$ of subspaces of complementary dimensions with $U\le W$ or $V=U\oplus W$, and $G$ contains a  graph automorphism;
		\item [\rm (4)]$M\cap X$ is $\Omega_{7}(2)$ or $\Omega_{8}^+(2)$ and $X$ is $P\Omega_{7}(q)$ or $P\Omega_{8}^+(q)$, respectively, $q$ is prime and $q\equiv \pm3\pmod{8}$;
 		\item [\rm (5)]$X=P\Omega_{8}^+(q)$, $q$ is prime and $q\equiv \pm3\pmod {8}$, $G$ contains a triality automorphism of $X$ and $M\cap X$ is $2^3\cdot 2^6\cdot PSL_{3}(2)$;
		\item [\rm (6)]$X=PSL_{2}(q)$ and $M\cap X$ is dihedral, $A_{4},S_{4},A_{5}$ or $PGL_{2}(q^\frac{1}{2})$;
		\item [\rm (7)]$X=PSU_{3}(5)$ and $M\cap X=M_{10}.$
	     \end{enumerate}
     	\end{enumerate}
	\end{enumerate}
\end{lemma}
\begin{lemma}\label{l6}{\rm (\cite[Theorem] {odd})}
Let $G$ be a primitive permutation group of odd degree $v$ on a set $\Omega$ with socle $X$, and $M=G_{\alpha}, \alpha \in \Omega$. If $X \cong A_m$, an alternating group, then one of the following holds:
\begin{enumerate}
	\item [\rm (1)]$M$ is intransitive, and $M=\left(S_s \times S_{t}\right) \cap G$ where $m=s+t,s\neq t$;
	\item [\rm (2)]$M$ is transitive and imprimitive, and $M=\left(S_s\wr S_t\right) \cap G$ where $m=st, s>1, t>1$;
	\item [\rm (3)]$M$ is primitive, $v=15$ and $G \cong A_7$ or $A_{8}$.
\end{enumerate}
\end{lemma}

\indent Next, we finish this section with some useful inequalities which will be used throughout the paper.
\begin{lemma}\label{l7}{\rm (\cite[Lemma 4.4] {large})}
	Suppose that $t$ is a positive integer. We have
	\begin{enumerate}
		\item [\rm (1)] if $t \ge 4$, then $t !<2^{4 t(t-3)/3}$;
		\item [\rm (2)]if $t \ge 5$, then $t !<5^{(t^2-3 t+1) / 3}$.
	\end{enumerate} 
\end{lemma}
\begin{lemma}\label{l8}{\rm (\cite[Lemma 4.2 and Corollary 4.3] {large})} 
	Let $G$ be a finite simple classical group over $\mathbb{F}_{q}$, with natural $($projective$)$ module $V$ of dimension $m$. Let $q=p^f$, where $p$ is a prime and $f\ge 1$. There are the following facts:
	\begin{enumerate}
		\item [\rm (1)]If $m \geqslant 2$, then 
		\begin{equation*}
			q^{m^2-2}<|PSL_{m}(q)|\le |SL_{m}(q)|<(1-q^{-2})q^{m^2-1},
		\end{equation*}
	\begin{equation*}
			(1-q^{-1})q^{m^2-2}<|PSU_{m}(q)|\le |SU_{m}(q)|<(1-q^{-2})(1+q^{-3})q^{m^2-1}.
		\end{equation*}
		\item [\rm (2)]If $m \geqslant 4$, then
	\begin{equation*}
		\frac{1}{2\beta}q^{m(m+1)/2}<|PSp_{m}(q)|\le |Sp_{m}(q)|\le (1-q^{-2})(1-q^{-4})q^{m(m+1)/2},
	\end{equation*}
with $\beta=(2,q-1)$.
		\item [\rm (3)]	If $m \geqslant 6$, then
		\begin{equation*}
			\frac{1}{8}q^{m(m-1)/2}<|P\Omega_{m}^\pm(q)|<|SO_{m}^\pm (q)|\le \delta (1-q^{-2})(1-q^{-4})(1+q^{-m/2})q^{m(m-1)/2},
		\end{equation*}
with $\delta =(2,q).$
	\end{enumerate}
\end{lemma}

\section{Proof of Theorem \ref{th1}}
In this section, let $\mathcal{D}$ be a 2-$(v,k,\lambda)$ non-trivial symmetric design of prime order, and $G\le Aut(\mathcal{D})$ be a flag-transitive almost simple group with scole $X$ is a finite non-abelian simple group. Note that $G$ is point-primitive, then the point stabilizer $M:=G_{\alpha}$ is maximal in $G$, and so $v=\frac{|X|}{|M\cap X|}$ from Lemma \ref{l3}. By \cite{class}, the finite non-abelian simple groups fall into four classes: (1) the 26 sporadic groups; (2) the alternating groups $A_{m}(m\ge 5)$; (3) the exceptional groups of Lie type; (4) the finite classical groups. \vspace{2mm}\\
{\bfseries Hypothesis 1.} Let $\mathcal{D}=(\mathcal{P}, \mathcal{B})$ be a non-trivial symmetric $2$-$(v,k,\lambda)$ design of prime order $n$, satisfying $(k,\lambda)=n$. Let $G$ be a flag-transitive automorphism group of $\mathcal{D}$ 
of almost simple type with socle $X$.\\
\indent In \cite{sporadic}, Tian and Zhou conclude that $X$ cannot be a sporadic simple group. Moreover, it is proved in \cite{exceptional} that $X$ also cannot be an exceptional group of Lie type. Hence, we have
\begin{proposition}\label{prop1}
	If $\mathcal{D}$ and $G$ satisfy Hypothesis $1$, then $X$ cannot be a sporadic simple group and finite exceptional group of Lie type.
\end{proposition}
 
 \indent Next, we mainly focus on the socle $X$ is an alternating group or a classical group.
\subsection{Classical groups}
In the following, assume that $X=X(q)$, the socle of $G$, is a non-abelian simple classical group, where $q=p^f$ for some positive integer $f$ and $p$ is a prime. Then $X=Soc(G)$ is isomorphic to one of the following groups:
\begin{enumerate}
	\item [\rm (1)]$PSL_{m}(q)$, with $m\ge 2$ and $(m,q)\neq (2,2),(2,3);$
	\item [\rm (2)]$PSU_{m}(q)$, with $m\ge 3$ and $(m,q)\neq (3,2);$
	\item [\rm (3)]$PSp_{2m}(q)$, with $m\ge 2$ and $(m,q)\neq (2,2);$
	\item [\rm (4)]$P\Omega_{2m+1}(q)$, with $m\ge 3$ and $q$ odd;
	\item [\rm (5)]$P\Omega_{2m}^\pm(q)$, with $m\ge 4$ .
\end{enumerate}

\indent  If $X=PSL_{2}(q)$, then $\mathcal{D}$ either is a $2$-$(7,4,2)$ design with $G=PSL_{2}(7)$ or $2$-$(11,6,3)$ design with $G=PSL_{2}(11)$ in \cite{psl}. Moreover, in \cite{four} we know that if $X=PSL_{m}(q)$ with $2<m\le 4$, then $\mathcal{D}$ is $2$-$(n^2+n+1,n^2,n^2-n)$ design, that is, the complement of $PG_{2}(n),$ and $PSL_{3}(n)\unlhd G\le P\Gamma L_{3}(n)$. 
  Now we can assume that the dimension $m$ of the underlying vector space is at least $5$. Recall that the order $n=k-\lambda=n(k^*-\lambda^*)$, so $k^*-\lambda^*=1$, and by \cite[Theorem]{lprime}, we may assume that $k^*=\lambda^*+1\ge 3.$ Furthermore, if $n=2$, then $k=2k^*$ and it follows from Lemma \ref{l1}(2)(i) that $k^*=2<3$, contracts the assumption. Hence, assume that $n$ is an odd prime divisor of $k.$ In what follows, we denote by $\hat{} M$ the pre-image of the group $M$ in the corresponding group, and the orders of outer automorphism groups are listed in Table \ref{t1}.
\begin{table}[!t]
	\centering
	\caption{Classical groups}\vspace{1mm}
	\begin{tabular}{l l l r}
		\hline $X=Soc(G)$ & $|X|$  & $d$& $|Ou t(X)|$ \\
		\hline$PSL_m(q)$ & $\frac{1}{d} q^{m(m-1) / 2} \prod_{i=2}^m\left(q^i-1\right)$ 
		 & $(m, q-1)$ &
		 \makecell [r]{$df \text { if } m=2$\\	$2df \text{ if } m \ge 3$}\\
		$PSU_m(q), m \geq 3$ & $\frac{1}{d} q^{m(m-1) / 2} \prod_{i=2}^m\left(q^i-(-1)^i\right)$ & $(m, q+1)$ & $ 2 d f$ \\
		$PSp_{2m}(q),m\ge 2$ & $\frac{1}{d} q^{m^2} \prod_{i=1}^m\left(q^{2 i}-1\right)$ & $(2, q-1)$ &  \makecell [r]{$(2, p)df \text { if } m=2$\\$df\text { if } m \geq 3 $
			} \\
		$P \Omega_{2 m+1}(q), m \geq 3$ & $\frac{1}{2} q^{m^2} \prod_{i=1}^m\left(q^{2 i}-1\right)$ & 2 & $2 f$ \\
		$P \Omega_{2 m}^{+}(q), m \geq 4$ & $\frac{1}{d} q^{m(m-1)}\left(q^m-1\right) \prod_{i=1}^{m-1}\left(q^{2 i}-1\right)$ & $\left(4, q^m-1\right)$ & \makecell[r]{	$6df \text { if } m=4$\\$2 df \text { if } m \geq 5$ 
		} \\
		$P \Omega_{2 m}^{-}(q), m \geq 4$ & $\frac{1}{d} q^{m(m-1)}\left(q^m+1\right) \prod_{i=1}^{m-1}\left(q^{2 i}-1\right)$ & $\left(4, q^m+1\right)$ & $ 2 d f$ \\
		\hline
	\end{tabular}
	\label{t1}
\end{table}
\begin{proposition}\label{prop2}
 If $\mathcal{D}$ and $G$ satisfy Hypothesis $1$, then $n$ divides $|M\cap X|_{p'}$. In particular, $|X|<2|M\cap X|\cdot |M\cap X|_{p'}^2$ and $|X|<|M\cap X|^3$. 
\end{proposition}
{\bfseries Proof}. Let $P(X)$ be the minimal degrees of permutation representations of finite classical simple groups $X$ in \cite{minimal}.\\
\indent From Lemmas \ref{l1}(3) and \ref{l3}, we know that $k\mid |Out (X)|\cdot|M\cap X|.$ Since $n$ is a prime divisor of $k$, then $n\mid |Out (X)|\cdot |M\cap X|.$ Suppose that $n\mid |Out (X)|\cdot |M\cap X|_{p}.$ We now consider each possibility for $X$ separately.\\
{\bfseries Case (1)}. $X=PSL_m(q)$, where $m\ge 5$ and $q=p^f$.\\
\indent By the main result of \cite{minimal},  $P(PSL_m(q))=\frac{q^m-1}{q-1}$, and then $v\ge \frac{q^m-1}{q-1}$. Since $n$ is an odd prime divisor of $|Out(X)|\cdot |X|_{p}$, by inspection $|Out(X)|$  and $|X|_{p}$ from Table \ref{t1}, we have $n\mid 2fp\cdot (m,q-1)$. So $n\le \text{max}\{f,p,q-1\}$, which implies that $n\le q$. From Lemma \ref{l1}(2)(i), we obtain that $k^*\le n$. Note that $v<2k$ by Lemma \ref{l1}(2)(ii), which follows that
\begin{equation*}
	\frac{q^m-1}{q-1}\le v<2k=2k^*n\le 2n^2\le 2q^2.
\end{equation*}
Recall that $m\ge 5$, the above inequality cannot hold.\\
{\bfseries Case (2)}. $X=PSU_m(q)$, where $m\ge 5$ and $q=p^f$.\\
\indent Firstly, assume that $(m,q)=(6s,2)$, where $s$ is an integer. Since $n\mid |Out(X)|\cdot |X|_{p}$, we have $n=3$. Combining the minimal degree $P(X)$ with Lemma \ref{l1}(2)(ii), so that $2^{6s-1}\cdot \frac{2^{6s}-1}{3}\le v<2k\le 18,$ which is impossible. Next, we deal with the case $(m,q)\neq (6s,2)$. Then $n\mid 2fp\cdot (m,q+1)$. Since $n$ is an odd prime, we conclude that $n\le \text{max}\{f,p,q+1\}\le q+1.$ Now, we have
\begin{equation*}
	\dfrac{(q^m-(-1)^m)(q^{m-1}-(-1)^{m-1})}{q^2-1}\le v<2k\le 2(q+1)^2,
\end{equation*}
and so $(q^5-1)(q^{4}-1)\le (q^m-(-1)^m)(q^{m-1}-(-1)^{m-1})<2(q+1)^2(q^2-1),$ which is impossible.\\
{\bfseries Case (3)}. $X=PSp_{2m}(q)$, where $m\ge 2$ and $q=p^f$.\\
\indent Assume first that $q>2$ and $(m,q)\neq (2,3)$. Since $n\mid |Out(X)|\cdot |X|_{p}$, we have $n\le \text {max}\{f,p\}\le q$. By the result of \cite{minimal}, we get $v\ge P(X)=\frac{q^{2m}-1}{q-1}$. Lemma \ref{l1}(2)(ii) implies that
\begin{equation*}
	\dfrac{q^{2m}-1}{q-1}\le v<2k<2q^2.
\end{equation*}
Note that $m\ge 2$, then $q^3+q^2+q+1<2q^2$, which is impossible. Next, if $m\ge 3$ and $q=2$, then $n=2$, a contradiction. If $(m,q)=(2,3),$ then $27\le v<2\cdot 3^2,$ which is also a contradiction.\\
{\bfseries Case (4)}. $X=P\Omega_{2m+1}(q)$, where $m\ge 3$ and $q=p^f$ odd.\\
\indent Assume that $q=3,$ then $f=1$. Since $n$ is an odd prime number dividing $|Out(X)|\cdot |X|_{p}$, then $n=3$. By Lemma \ref{l1}(2)(ii), $v<2k$, we have $\frac{3^m(3^m-1)}{2}\le v<2\cdot 3^2$,
 which is impossible for $m\ge 3$. If $q\ge 5,$ then $v\ge P(X)=\frac{q^{2m}-1}{q-1}$ from \cite{minimal}. By inspection $|Out(X)|$ and $|X|_{p}$, we have $n\mid 2fp$, and so $n\le q$. From Lemma \ref{l1}(2), we obtain that $\frac{q^{2m}-1}{q-1}\le v<2k\le 2q^2,$ contrary to $m\ge 3$.\\
{\bfseries Case (5)}. $X=P\Omega_{2m}^\epsilon(q)$ with $m\ge 4$ and $q=p^f$, where $\epsilon=\{+,-\}$. \\
 \indent Suppose that $X=P\Omega_{2m}^+(2).$ Then $v\ge 2^{m-1}(2^m-1)$, the minimal degree of $X$, by \cite{minimal}. Since $n$ is an odd prime divisor of $|Out(X)|\cdot |X|_{p}$, then we have $n=3.$ Now Lemma \ref{l1}(2)(ii) says that $v<2k$, that is, $ 2^{m-1}(2^m-1)<2\cdot 3^2$, which is impossible since $m\ge 4$. When $X\neq P\Omega_{2m}^+(2)$, we know the minimal degree $P(X)\ge \frac{(q^m+1)(q^{m-1}-1)}{q-1}$. Since $n$ is an odd prime number and $n\mid |Out(X)|\cdot |X|_{p}$, then $n\le q.$ Combining it with Lemma \ref{l1}(2), we conclude that
 \begin{equation*}
 \dfrac{(q^m+1)(q^{m-1}-1)}{q-1}\le v <2k\le 2q^2.
 \end{equation*}
Recall that $m\ge 4$, then $(q^4+1)(q^3-1)<2q^2(q-1),$
 a contradiction.\\
 \indent Hence, $(n,|Out(X)|\cdot |M\cap X|_{p})=1$ and $n\mid |M\cap X|_{p'}$. Recall that $G$ is point-primitive, then $v=|X:M\cap X|$
 by Lemma \ref{l3}. Furthermore, Lemma \ref{l1}(2) imply that $|X:M\cap X|=v<2k\le 2n^2\le 2|M\cap X|_{p'}^2.$ By simple computation, we have $|X|< 2|M\cap X||M\cap X|_{p'}^2\le |M\cap X|^3.\hfill \square$\vspace{3mm}\\
 \indent Recall that Lemma \ref{l2}, we know $G$ is a primitive permutation group of odd degree on point set. Next, we will apply Proposition \ref{prop2} to analyze all cases in Lemma \ref{l5}. Proposition \ref{prop2} is very useful for proof of the main theorem.
 \begin{proposition}\label{prop3}
 	If $\mathcal{D}$ and $G$ satisfy Hypothesis $1$, then $X\neq PSL_{m}(q)$ for $m\ge 5$.
 \end{proposition}
{\bfseries Proof}. Putting $M_0=M \cap X$, where $M=G_\alpha$ with $\alpha \in \mathcal{P}$. In the following, we analyse each of these possible cases in Lemma \ref{l5} respectively.\vspace{2mm}\\
{\bfseries Case (1)}. $M_0$ is a parabolic subgroup of $X$.\\
\indent In this case, $M=P_{i}$ stabilises a subspace of $V$ of dimension $i$, where $m\ge 2i$. Suppose first that $M=P_{1}$, then $G$ is $2$-transitive. Recall that all $2$-transitive symmetric designs have been determined by Kantor \cite{kantor}. By \cite[Theorem]{kantor}, we know that $\mathcal{D}$ and $G$ is one of the following: (i) the complement of the unique Hadamard design $2$-$(11,5,2)$, and $G=PSL_{2}(11)$; (ii) the complement of a projective plane $PG_{2}(n)$, and $PSL_{3}(n)\unlhd G\le P\Gamma L_{3}(n)$. But none of the above parameters are satisfied $m\ge 5.$\\
\indent Now suppose $M=P_{i}$ with $1<i<m$ fixed an $i$-subspace of $V$. By \cite[Proposition 4.1.17]{classical}, we know that $M_{0}\cong \hat{}[q^{i(m-i)}]:SL_{i}(q)\times  SL_{m-i}(q)\cdot (q-1).$ Combining it with Proposition \ref{prop2}, then $n\le \frac{q^{m-i}-1}{q-1}$. On the other hand, we have $v>q^{i(m-i)}$. Moreover, Lemma \ref{l1}(2) follows that
\begin{equation}\label{e3.3.1}
q^{i(m-i)}<v<2k\le 2n^2\le \dfrac{2(q^{m-i}-1)^2}{(q-1)^2}.\tag{3.3.1}
\end{equation}
Note that $\frac{q^{j}-1}{q-1}<2q^{j-1}$ for some integer $j$, then (\ref{e3.3.1}) implies that $q^{i(m-i)}<q^{2m-2i+1}$, and hence $m(i-2)<i^2-2i+1$. Recall that $m\ge 2i$, we have $2i(i-2)<i^2-2i+1,$ and so $i^2<2i+1.$ Thus $i=2$. If $i=2$, then $G$ is a primitive rank 3 group, with non-trivial subdegrees $d_{1}=\frac{q(q+1)(q^{m-2}-1)}{q-1}$ and $d_{2}=\frac{q^4(q^{m-2}-1)(q^{m-3}-1)}{(q^2-1)(q-1)}$ by \cite{linear}. The symmetric 2-designs admitting primitive rank 3 automorphism groups have been classified by Dempwolff \cite[Theorem]{De}. It is easy to observe that there are no parameters that satisfy $n=k-\lambda$ is a prime. \vspace{2mm}\\
{\bfseries Case (2)}. $M$ is preserve a decomposition  $V=V_1 \oplus \cdots \oplus V_t$  with each $V_{j}$ of the same dimension $i$, where $m=it$. In addition, $q$ is an odd number.\\
\indent By \cite[Proposition 4.2.9]{classical}, we have $M_{0}\cong \hat{}  SL_{i}(q)^t \cdot (q-1)^{t-1} \cdot S_{t}.$ Furthermore, we conclude that $v>\frac{q^{m(m-i)}}{t!}$. If $i=1$, then $m=t$. In this case, we know $|M\cap X|=\frac{(q-1)^{m-1}\cdot m!}{(m,q-1)}.$ Then Lemma \ref{l8} and Proposition \ref{prop2} imply that 
\begin{equation}\label{e3.3.2}
q^{m^2-2}<|X|<|M\cap X|^3=\frac{(m!)^3(q-1)^{3m-3}}{(m,q-1)^3}\tag{3.3.2}
\end{equation}
Since $q$ is odd, then $3^{m^2-3m+1}\le q^{m^2-3m+1}\le (m!)^3$ by (\ref{e3.3.2}). It follows that $m=5$, and again from (\ref{e3.3.2}), we have $q^{11}<2^9\cdot 3^3\cdot 5^3$. Hence, $q=3$. From (\ref{e3.3.2}), we obtain that $3^{23}<(5!)^{3}\cdot 2^{12},$ a contradiction. Thus $i\ge 2.$\\
\indent Since $n$ is an odd prime divisor of $|M_{0}|_{p'}$, then $n\le \text{max} \{t,\frac{q^i-1}{q-1}\}$, and so $n\le \frac{t(q^i-1)}{q-1}$. Lemma \ref{l1}(2) follows that
	$\frac{q^{m(m-i)}}{t!}<v<2k\le 2n^2\le \frac{2t^2(q^i-1)^2}{(q-1)^2},$
and
\begin{equation}\label{e3.3.3}
	q^{m(m-i)}<2^3\cdot t^2\cdot q^{2i-2}\cdot t!.\tag{3.3.3}
\end{equation}
 Assume that $t=2$. Then $m=2i$ and $q^{2i^2}<2^6\cdot q^{2i-2}$ by (\ref{e3.3.3}). It implies that $i^2<i+2$, a contradiction. If $t=3$, then $m=3i$. Similarly, by (\ref{e3.3.3}), we conclude that $q^{6i^2}<2^4\cdot 3^3\cdot q^{2i-2},$ and so $q^{6i^2-2i+2}<2^4\cdot 3^3<3^6$. Note that $q$ is odd, then $6i^2-2i+2<6$, and so 
 $3i^2<i+2,$ a contradiction. If $t\ge 4,$ by Lemma \ref{l7}(1), $t !<2^{4 t(t-3)/3}$. From 
 (\ref{e3.3.3}), we have $q^{3m(m-i)}<2^{4t^2-12t+9}\cdot t^6\cdot q^{6i-6}$. Since $t^6\le 2^{4t}$, then $q^{3m(m-i)}<2^{4t^2-8t+9}\cdot q^{6i-6}$. Combining it with $m=it$ implies  that $t^2(3i^2-4)+t(8-3i^2)<6i+3,$ and then $t(3i^2-4)+8-3i^2<6i$. Recall that $t\ge 4,$ we get $9i^2<6i+8$, which is impossible for $i\ge2.$\vspace{2mm}\\
{\bfseries Case (3)}. Assume that $M$ be stabilize a pair $\{U,W\}$ of subspaces $U$ and $W$, where $U$ has dimension $i$ and $W$ has dimension $m-i$ with $1\le i<\frac{m}{2}$. Moreover, either $U \leqslant W$ or $U \oplus W=V$. In addition, $G$ contains a graph automorphism.\\
\indent Firstly, suppose that $U \leqslant W$, then $M_{0} \cong \hat{ }$$[q^{2im-3i^2}]\cdot SL_{i}(q)^2 \times SL_{m-2i}(q) \cdot (q-1)^2$ by \cite[Proposition 4.1.22]{classical}. Since \begin{equation*}
	v=\dfrac{q^{m(m-1)/2}\cdot (q^{i+1}-1)\cdots (q^{m-i}-1)(q^{m-i+1}-1)\cdots (q^{m}-1)}{q^{2im-3i^2}\cdot q^{i(i-1)}\cdot q^{(m-2i)(m-2i-1)/2}\cdot (q-1)\cdots (q^{m-2i}-1)(q-1)\cdots (q^i-1)},
\end{equation*}
 we have $v>q^{i(2m-3i)}.$ There is a subdegree which is a power of $p$ by Lemma \ref{l4}. Moreover, note that $q$ is odd, then $(v-1)_{p}=q$. Then by Lemma \ref{l1}(4), we have $k^*\mid q.$ Since $n$ is an odd prime divisor of $|M_{0}|_{p'},$ then $n\le q^{m-2}-1$. Hence, Lemma \ref{l1}(2)(ii) implies that 
\begin{equation}\label{e3.3.4}
	q^{i(2m-3i)}<v<2k=2k^*n\le 2q(q^{m-2}-1)<q^m,\tag{3.3.4}
\end{equation} 
and $m(2i-1)<3i^2$. Note that $m>2i$, then $2i(2i-1)<3i^2,$ and so $i=1.$ In this case, we have $q^{2m-3}<q^m$ from (\ref{e3.3.4}). Then $2m-3<m$, a contradiction.\\
\indent If $U \oplus W=V$, then $M_{0}\cong \hat{ }SL_{i}(q)\times SL_{m-i}(q)\cdot (q-1)$ by \cite[Proposition 4.1.4]{classical}. By simple calculation, we have $v>q^{2i(m-i)}.$
 On the other hand, by Proposition \ref{prop2}, we get $n\le q^{m-i}-1$. Combining this with $v>q^{2i(m-i)}$, we obtain that $q^{2i(m-i)}<v<2k\le 2n^2< 2q^{2m-2i}.$ It follows that $m(2i-2)<2i^2-2i+1$. Note that $m>2i$, then $ 2i(2i-2)<m(2i-2)<2i^2-2i+1$, and hence $2i^2<2i+1.$ Thus $i=1$, and $v=q^{m-1}(q^m-1)/(q-1)$. Since $n$ is an odd prime divisor of $|M_{0}|_{p'}$, then $n$ divides $q^j-1$ with $j\le m-1$. Hence, $n\le \frac{q^{m-1}-1}{q-1}$. Furthermore, 
\begin{equation*}
	\dfrac{q^{m-1}(q^m-1)}{q-1}=v<2n^2\le 2\dfrac{(q^{m-1}-1)^2}{(q-1)^2},
\end{equation*}
and so $q^{m-1}(q^m-1)(q-1)<2(q^{m-1}-1)^2<2\cdot q^{2m-2}.$ It follows that $(q^m-1)(q-1)<2q^{m-1},$ which is impossible since $q$ is odd.\vspace{2mm}\\
{\bfseries Case (4)}. Assume that $q=q_0^t$ is odd with $t$ odd prime, and $M=N_G(X)$.\\ \indent According to \cite[Proposition 4.5.3]{classical}, $M\cap X\cong (m,\frac{q-1}{q_{0}-1})\cdot \hat{ } SL_{m}(q_{0}).$ By Lemma \ref{l8}, we have $v>q_{0}^{t(m^2-2)-m^2+1}$. On the other hand, Proposition \ref{prop2} says that $n\le q_{0}^m-1$. It follows from that 
\begin{equation*}
	q_{0}^{t(m^2-2)-m^2+1}<v<2n^2\le 2(q_{0}^m-1)^2.
\end{equation*} 
 and so $t(m^2-2)-m^2+1<2m+1.$ As $t$  is odd prime, it is easy to know that $m^2<m+3,$ a contradiction.$\hfill \square$ \vspace{3mm} 
\begin{proposition}\label{prop4}
	If $\mathcal{D}$ and $G$ satisfy Hypothesis $1$, then 
  $X\neq PSU_{m}(q)$ for $m\ge 5.$
\end{proposition}
{\bfseries Proof}. Let $M_0=M \cap X$, where $M=G_\alpha$ with $\alpha \in \mathcal{P}$. 
We analyse each of these cases in Lemma \ref{l5} separately.\vspace{2mm}\\
{\bfseries Case (1)}. Assume that $M\cap X$ is the parabolic subgroup of $X$ with $q=2^f$.\\
\indent From \cite[Proposition 4.1.18]{classical}, we know $M_{0}\cong \hat{ } q^{i(2m-3i)}:SL_{i}(q^2)\cdot SU_{m-2i}(q)\cdot (q^2-1),$ where $i\le [m/2]$. According to \cite[Section 6]{linear}, we have $v>q^{i(2m-3i)}.$ On the other hand, there is a subdegree which is a power of $2$ by Lemma \ref{l4}. Moreover, the $2$-part of $v-1$ is one of the following possibilities:
\begin{enumerate}
	\item [\rm (a)]if $m$ is odd and $i=\frac{m-1}{2}$, then $(v-1)_{2}=q^3;$
	\item [\rm (b)]if $m$ is even and $i=\frac{m}{2}$, then $(v-1)_{2}=q;$
	\item [\rm (c)]otherwise,  $(v-1)_{2}=q^{2}.$
\end{enumerate}
Then $k^*\mid q^3$ by Lemma \ref{l1}(4). Proposition \ref{prop2} says that $n\le q^{m-2}+1$. Recall that $v<2k$, then 
\begin{equation}\label{e3.4.1}
	q^{i(2m-3i)}<2q^3(q^{m-2}+1)<2^2\cdot q^{m+1},\tag{3.4.1}
\end{equation}
and so $m(2i-1)<3i^2+3.$ Note that $m\ge 2i,$ we have $2i(2i-1)<3i^2+3.$ It implies that $i^2<2i+3$, and hence $i=1$ or $2$. If $i=1$, then $(v-1)_{2}=q^2$ by the condition $m\ge 5$. Then by (\ref{e3.4.1}), we get $q^{2m-3}<2q^2(q^{m-2}+1)<2^2\cdot q^m$. Since the value of $m$ is at least $5$, then $2^2\le q^2\le q^{m-3}<2^2,$ a contradiction. Hence $i=2$, by (\ref{e3.4.1}), we have $q^{4m-12}<2^2\cdot q^{m+1}$. It implies that $q^{3m-13}<2^2$, which is impossible for $m\ge 5$.\vspace{2mm}\\
{\bfseries Case (2)}. Suppose that $M$ is the stabilizer of a non-singular subspace and $q$ is odd.\\
\indent By \cite[Proposition 4.1.4]{classical},  $M_{0}\cong \hat{ } SU_{i}(q)\times SU_{m-i}(q)\cdot (q+1)$, where $m>2i$. By Lemma \ref{l8}, we have $v>q^{2i(m-i)-6}$. Note that $n$ is an odd prime, Proposition \ref{prop2} says that $n\mid |M_{0}|_{p'}$, then $n\le q^{m-i}+1$. By Lemma \ref{l1}(2)(ii), $v<2k$, we conclude that $q^{2i(m-i)-6}<2(q^{m-i}+1)^2.$ Note that $q^j+1<\frac{3}{2}q^j$ for some integer $j$, then $q^{2i(m-i)-6}<q^{2m-2i+2}$, and so $m(2i-2)<2i^2-2i+8.$ As $m>2i$, we obtain that $i^2<i+4.$ Therefore, $i=1$ or $2.$\\
\indent Assume that $i=1$, in which case we know $v=\frac{q^{m-1}(q^{m}-(-1)^m)}{q+1}.$ If $m$ is even, then $v=\frac{q^{m-1}(q^{m}-1)}{q+1}=q^{m-1}(q^{m-1}-q^{m-2}+\cdots +q-1)$ is even, a contradiction. Thus $m$ is odd, and $v=\frac{q^{m-1}(q^{m}+1)}{q+1}$. By Lemma \ref{l1}(2), we have  $\frac{q^{m-1}(q^{m}+1)}{q+1}<2\cdot \frac{(q^{m-1}-1)^2}{(q-1)^2}$, and so 
\begin{equation*}
	q^{2m-1}(q-1)^2<q^{m-1}(q^{m}+1)(q-1)^2<2(q^{m-1}-1)^2(q+1)<2q^{2m-2}(q+1).
\end{equation*}
It follows that $q(q-1)^2<2(q+1),$ which is impossible. Next, let $i=2$. Then $v={q^{2m-4}(q^{m}-(-1)^m)(q^{m-1}-(-1)^{m-1})}/{(q^2-1)(q+1)}$. In the same way, we have
\begin{equation*}
	\dfrac{q^{2m-4}(q^{m}-(-1)^m)(q^{m-1}-(-1)^{m-1})}{(q^2-1)(q+1)}<2\dfrac{(q^{m-2}-1)^2}{(q-1)^2}<\dfrac{2q^{2m-4}}{(q-1)^2},
\end{equation*}
and so $(q^5+1)(q^4-1)(q-1)\le (q^{m}-(-1)^m)(q^{m-1}-(-1)^{m-1})(q-1)<2(q+1)^2$, a contradiction.\vspace{2mm}\\
{\bfseries Case (3)}. Assume that $q$ is odd, and $M \cap X$ is the stabilizer of an orthogonal decomposition $V=V_1 \oplus \cdots \oplus V_t$ with $\operatorname{dim}\left(V_j\right)=i$, where $m=it$ $(t\ge 2)$.\\
\indent By \cite[Proposition 4.2.9]{classical}, we know $M_{0}$ is isomorphic to $\hat{ } SU_{i}(q)^t\cdot (q+1)^{t-1}\cdot S_{t}.$ Firstly, we show that $i\neq 1.$ In fact, if $i=1$, then $M_{0}\cong \hat{ }(q+1)^{m-1}\cdot S_{m}.$ From Lemma \ref{l1}(2) and Proposition \ref{prop2}, we obtain that 
\begin{equation*}
	\dfrac{q^{\frac{m(m-1)}{2}}\cdot \prod_{j=2}^m(q^j-(-1)^j)}{(q+1)^{m-1}\cdot m!}=v<2k\le 2n^2\le \frac{m^2(q+1)^2}{2}.
\end{equation*}
Note that $q^{\frac{m(m-1)}{2}}<\prod_{j=2}^m(q^j-(-1)^j)$,
we conclude that $q^{m(m-1)}<2^m\cdot q^{m+1}\cdot m^2\cdot m!.$ Since $m\ge 5$, by Lemma \ref{l7}(2), then $q^{m(m-1)}<2^m\cdot q^{m+1}\cdot m^2\cdot 5^{\frac{m^2-3m+1}{3}}.$ Combining it with $m^6<3^{2m},$ yields that $q^{3m^2-3m}<q^{2m^2+m+5}.$ Then $3m^2-3m<2m^2+m+5$, and so $m^2<4m+5$, which is impossible for $m\ge 5.$ Thus, $i\ge 2$.\\
\indent If $i\ge 2,$ then $v>{q^{{i^2t(t-1)}/{2}}}/{t!}$. Proposition \ref{prop2} implies that $n\le t\cdot (q^i-(-1)^i).$ Similarly, we have 
\begin{equation}\label{e3.4.2}
	\frac{q^{i^2t(t-1)/2}}{t!}<v\le 2n^2\le 2t^2(q^i-(-1)^i)^2.\tag{3.4.2}
\end{equation}
If $t=2$, then $n\le q^i-(-1)^i$, and so $q^{i^2}<2^6\cdot q^{2i}<q^{2i+4}.$ It implies that $i^2<2i+4$, and then $i=2$ or $3$. If $i=2$, then $m=it=4<5$, a contradiction.  Hence $i=3$, by (\ref{e3.4.2}), we have $q^9<2^6\cdot q^6$, and then $q=3.$ Now again from (\ref{e3.4.2}), we get $3^9<2^8\cdot 7^2, $ a contradiction. Assume that $t=3$, then $n\le q^i-(-1)^i$ from Proposition \ref{prop2}. By the inequality (\ref{e3.4.2}), we obtain that $q^{3i^2}<q^{2i+5}$, and so $3i^2<2i+5$.
It follows that $i=1$, which is a contradiction and thus $t\ge 4.$ Now, by Lemma \ref{l7}(1) and (\ref{e3.4.2}), we have
\begin{equation*}
	q^\frac{t(t-1)i^2}{2}<2^3\cdot t^2\cdot 2^\frac{4t(t-3)}{3}\cdot q^{2i} .
\end{equation*}
Note that $t^3\le 3^t,$ then $q^{3t(t-1)i^2}<2^{8t(t-3)+18}\cdot t^{12}\cdot q^{12i}<q^{12i+6t^2-14t+12}.$ It follows that $t^2(3i^2-6)<t(3i^2-14)+12i+12<t(3i^2-14)+12it.$ So $t(3i^2-6)<3i^2-14+12i$. Moreover, since $t\ge 4,$
 then $9i^2<12i+10$, a contradiction.\vspace{2mm}\\
{\bfseries Case (4)}. Suppose that $q=q_0^t$ is odd with $t$ odd prime, and $M=N_G(X)$. \\
\indent By \cite[Proposition 4.5.3]{classical}, we know that $M_{0}\cong \hat{ }SU_{m}(q_{0})\cdot (\frac{q+1}{q_{0}+1},m)$. Here by Proposition \ref{prop2}, we conclude that $|X|<|M_{0}|^3,$ that is,
\begin{equation*}
	q_{0}^{tm(m-1)}<(q+1)\cdot q_{0}^{3t-3}\cdot q_{0}^{3m^2-3}<2q_{0}^{3m^2+4t-6},
\end{equation*}
which implies that $tm(m-1)<3m^2+4t-5$, and hence $t(m^2-m-4)<3m^2-5$. It follows from the fact that $t$ is odd prime, we have $t=3$. Then $q=q_0^3$, and Lemma \ref{l8} implies that
$v>q_{0}^{2m^2-10}.$ On the other hand, by Proposition \ref{prop2}, we know $n$ is an odd prime divisor of $|M_{0}|_{p'}$, and so $n\le q_{0}^m+1.$ Now Lemma \ref{l1}(2) says that 
\begin{equation*}
	q_{0}^{2m^2-10}<v<2n^2\le 2\cdot (q_{0}^m+1)^2<2^3\cdot q_{0}^{2m}<q_{0}^{2m+2}.
\end{equation*}
Thus $2m^2-10<2m+2$, that is, $m^2<m+6,$
which is impossible for $m\ge 5.$$\hfill \square$ \vspace{3mm} 
\begin{proposition}\label{prop5}
	If $\mathcal{D}$ and $G$ satisfy Hypothesis $1$, then 
$X$ cannot be $PSp_{2m}(q)$ with $m\ge 2$ and $(m,q)\neq (2,2),(2,3).$
\end{proposition}
{\bfseries Proof}.  Let $M_0=M \cap X$, where $M=G_\alpha$ with $\alpha\in \mathcal{P}$. Similarly, we analyse each of the cases in Lemma \ref{l5}.\vspace{2mm}\\
{\bfseries Case (1)}. Assume that $M_{0}$ is the parabolic subgroup of $X$ with $q=2^f$.\\
\indent  In this case, there is a subdegree which is a power of $2$ by Lemma \ref{l4}. On the other hand, we know that $(v-1)_{2}=q.$  Then Lemma \ref{l1}(4) gives $k^*\mid q$. By \cite[Proposition 4.1.19]{classical}, we have $M_{0}\cong [q^h]\cdot (GL_{i}(q)\times PSp_{2m-2i}(q))$, where $h=2mi+\frac{i-3i^2}{2}$ and $i\le m.$ According to \cite[Proposition 4.3]{lprime}, we get $v>q^{i(4m-3i)}.$ By the fact that $n$ is an odd prime divisor of $|M_{0}|_{p'}$, then $n\le \frac{q^m-1}{q-1}$. Again by Lemma \ref{l1}(2), we have 
\begin{equation}\label{e3.5.1}
	q^{i(4m-3i)}<2k\le 2\cdot q\cdot \frac{q^m-1}{q-1}<2^2\cdot q^m\le q^{m+2},\tag{3.5.1}
\end{equation} 
it follows that $i(4m-3i)<m+2.$ Then $m(4i-1)<3i^2+2,$ and so $i^2<i+2$. The inequality holds only for $i=1$, in which case  $q^{4m-3}<q^{m+2}$ by (\ref{e3.5.1}), the desired contradiction.\vspace{2mm}\\
{\bfseries Case (2)}. Assume that $q$ is odd, and $M$ is the stabilizer of a non-singular subspace.\\
\indent  According to \cite[Proposition 4.1.3]{classical}, $M_{0}\cong  2.PSp_{2i}(q)\times PSp_{2m-2i}(q)$, where $m>2i$. Moreover, we have $v>q^{4i(m-i)}.$
 By Proposition \ref{prop2}, $n$ divides $|M_{0}|_{p'}$, and so does $n\le q^{m-i}+1$. Now $q^{4i(m-i)}<v\le 2n^2<2(q^{m-i}+1)^2,$ which implies that $2i(m-i)<m-i+1.$ As $m>2i$, then $2i(2i-1)<m(2i-1)<2i^2-i+1,$ that is,
$2i^2<i+1$, a contradiction. \vspace{2mm}\\
{\bfseries Case (3)}. Suppose that $q$ is odd, and $M_{0}$ is the stabilizer of an orthogonal decomposition $V=V_1 \oplus \cdots \oplus V_t$ with $\operatorname{dim}\left(V_j\right)=i$, where $m=it$ $(t\ge 2)$.\\
\indent Now, we know that $v>{q^{2i^2t(t-1)}}/{(t!)}.$ Proposition \ref{prop2} implies that $n\le t(q^i+1).$ We apply Lemma \ref{l1}(2) and conclude that
\begin{equation}\label{e3.5.2}
	q^{2i^2t(t-1)}<2t^2(q^i+1)^2\cdot t!\tag{3.5.2}
\end{equation}
 If $t=2$, then $q^{2i^2}<2^3\cdot q^{i}$ from (\ref{e3.5.2}). It follows that $2i^2<i+2$. The previous inequality holds only for $i=1$ in which case $v=\frac{q^2(q^2+1)}{2}$. From Proposition \ref{prop2}, we have $n\le \frac{q+1}{2}$. Then $\frac{q^2(q^2+1)}{2}\le 2\cdot \frac{(q+1)^2}{4}$ by Lemma \ref{l1}(2). It implies that $q^4<2q+1$, contrary to $q\ge 3$. If $t=3$, then $q^{12i^2}<432q^{2i},$ and so $6i^2<i+3.$ This is impossible. Finally, assume that $t\ge 4.$ By Lemma \ref{l7}(1), $t !<2^{4 t(t-3)/3}$, we have
 $q^{6i^2t(t-1)}<t^6\cdot 2^{4t(t-3)+9}\cdot q^{6i}.$ Note that $t^6\le 3^{2t}$, we have $3i^2t(t-1)<3+t+3i+2t(t-3)$, that is, $t^2(3i^2-2)<t(3i^2-5)+3i+3.$ It follows that $t^2(3i^2-2)<t(3i^2-5)+3it,$ and so $t(3i^2-2)<3i^2-5+3i.$
 Since $t\ge 4$, then $3i^2<i+1$, which is impossible.\vspace{2mm}\\
{\bfseries Case (4)}. Finally, let $q=q_0^t$ is odd with $t$ odd prime, and $M=N_G(X)$.\\
  \indent By \cite[Proposition 4.5.4]{classical}, we know that $M_{0}\cong PSp_{2m}(q_{0})$. From Proposition \ref{prop2}, $|X|<|M_{0}|^{3},$ we get
\begin{equation*}
	q^{m^2} \cdot q^{m^2}\le	q^{m^2} \cdot \prod_{i=1}^{m}(q^{2i}-1)<q_{0}^{3m^2}\cdot (\prod_{i=1}^{m}(q_{0}^{2i}-1))^3\le 
	q_{0}^{6m^2+3m}.
\end{equation*}
 It follows that $2tm^2<6m^2+3m,$ and hence $t=3$. If $t=3$, by Lemma \ref{l8}, we have
 $v>q_{0}^{4m^2+2m-2}$ by Lemma \ref{l8}.
  Since $n$ is an odd prime divisor of $|M_{0}|_{p'}$, then $n\le q_{0}^m+1.$ By Lemma \ref{l1}(2), we conclude that $q_{0}^{4m^2+2m-2}<2(q_{0}^m+1)^2.$ Note that $q_{0}^m+1<2q_{0}^m,$ then $q_{0}^{4m^2-2}<8,$ which is impossible for $m\ge 2.$$\hfill \square$ \vspace{3mm} 
\begin{proposition}\label{prop6}
	If $\mathcal{D}$ and $G$ satisfy Hypothesis $1$, then 
	 $X=Soc(G)$ cannot be $P\Omega_{2 m+1}(q)$ with $q$ odd.
\end{proposition}
{\bfseries Proof}. Let $M_0=G_\alpha \cap X$, where $\alpha\in \mathcal{P}$.  Note that $\Omega_{3}(q)\cong PSL_{2}(q)$ and  $\Omega_{5}(q)\cong PSp_{4}(q)$, then we may assume that $m\ge 3.$ If $M=N_{G}(X))$ with $q=q_{0}^t$ odd and $t$ odd prime, then we argue exactly the same as in symplectic groups.
 Except for this case, we analyse other cases in the following.\vspace{2mm}\\
{\bfseries Case (1)}. Assume that $M$ is the stabilizer of a non-singular subspace and $q$ is odd. In this case, $M=N_{i}^\epsilon$ with $m\ge i.$\\
\indent If $i=1$, then $M_{0}\cong \hat{ } \Omega_{2m}^{\epsilon}(q)\cdot 2,$ where $\epsilon=+$ or $-$. If $\epsilon=+$, then $v=\frac{q^m(q^m+1)}{2}$. By Proposition \ref{prop2}, we know that $n\le \frac{q^m-1}{q-1}$. Now by Lemma \ref{l1}(2), then $\frac{q^m(q^m+1)}{2}=v<2n^2<\frac{2(q^m-1)^2}{(q-1)^2}$.
It follows that $q^m(q^m+1)(q-1)^2<4(q^m-1)^2$, and so $(q-1)^2<4$, which is impossible for $q\ge 3.$ In the same way, if $\epsilon=-$, then $v=\frac{q^m(q^m-1)}{2}$. Moreover, we have
\begin{equation}\label{e3.6.1}
	\dfrac{q^m(q^m-1)}{2}=v<2n^2<\dfrac{2(q^m+1)^2}{(q+1)^2},\tag{3.6.1}
\end{equation}  
it follows that $q^m(q^m-1)(q+1)^2<4(q^m+1)^2.$ Note that
$\frac{1}{2}q^m<q^m-1<q^m+1<\frac{3}{2}q^m$, we get $(q+1)^2<18$, and thus $q=3$. If $q=3$, then $3^{2m+1}<6\cdot 3^m+1$ from (\ref{e3.6.1}). Then $2m+1<m+2$, a contradiction. Hence, $i\ge 2.$\\
\indent If $i\ge 2,$ by \cite[Proposition 4.1.6]{classical}, then $M_{0}\cong \Omega_{i}^\epsilon(q)\times \Omega_{s-i}(q)\cdot 4$, where $\epsilon\in \{+,-\}$, $i$ is even and $s=2m+1.$ Since $\frac{q^j-1}{q^\ell-1}>q^{j-\ell}$ for some integers $j$ and $\ell$, then $v>\frac{q^{i(s-i)}}{4}$. Hence by Proposition \ref{prop2}, we have  $n\le \frac{q^{(s-i-1)/2}+1}{2}.$ Again by Lemma \ref{l1}(2),
\begin{equation*}
	\dfrac{q^{i(s-i)}}{4}<v<2\cdot \dfrac{(q^{(s-i-1)/2}+1)^2}{4}=\dfrac{(q^{(s-i-1)/2}+1)^2}{2}.
\end{equation*}
It implies that $s(i-1)<i^2-i+1.$ Recall that $s=2m+1$ and $m\ge i$, then $s>2i.$ Thus $2i(i-1)<i^2-i+1,$ and $i^2<i+1,$ which is impossible.\vspace{2mm}\\
{\bfseries Case (2)}. Assume that $q$ is odd, and $M_{0}$ is the stabilizer of an orthogonal decomposition $V=V_1 \oplus \cdots \oplus V_t$ with $\operatorname{dim}\left(V_j\right)=i$, where $i$ is odd.\\
\indent Let $i=1.$ By \cite[Proposition 4.2.15]{classical}, we know $M_{0}\cong 2^{2m}\cdot S_{2m+1}$ if $q\equiv \pm 1 \pmod{8}$ or $2^{2m}\cdot A_{2m+1}$ if $q\equiv \pm3 \pmod{8}$, respectively. From Proposition \ref{prop2} , $|X|<2|M\cap X|\cdot |M\cap X|_{p'}$, we have the possible of $X$ and $M_{0}$ which are listed in Table \ref{t2}. Note that for each case the maximal value of $n$ is given in the last column of Table \ref{t2}. However, all of these possibilities have $v>2n^2\ge 2k$, a contradiction.\\
\begin{table}[!t]
	\centering
	\caption{Possibilities of $(X,M_{0})$ in Proposition \ref{prop6} Case (2) when $i=1$}\vspace{1mm}
	\begin{tabular}{c l l l l r}
		\hline Line & $X$ & $M_{0}$ & $|X|$&$v=|X:M_{0}|$ & $n\le$ \\
		\hline $1$ & $\Omega_{7}(3)$ & $2^6\cdot A_{7}$ & $2^9\cdot 3^9\cdot 5\cdot 7 \cdot 13$&$3^7\cdot 13$ & $7$\\
		 $2$ & $\Omega_{7}(5)$ & $2^6\cdot A_{7}$ & $2^9\cdot 3^4\cdot 5^9\cdot 7 \cdot 13\cdot 31$&$3^2\cdot 5^8\cdot 13\cdot 31$ & $7$\\
		 $3$ & $\Omega_{9}(3)$ & $2^8\cdot A_{9}$ & $2^{14}\cdot 3^{16}\cdot 5^2\cdot 7 \cdot 13\cdot 41$&$3^{12}\cdot 5\cdot 13\cdot 41$ & $7$\\
		 $4$ & $\Omega_{11}(3)$ & $2^{10}\cdot A_{11}$ & $2^{17}\cdot 3^{4}\cdot 5^2\cdot 7 \cdot 11$&$3^{21}\cdot 11\cdot 13\cdot 41\cdot 61$ & $11$\\
		 \hline
	\end{tabular}
	\label{t2}
\end{table}
\indent If $i=3$, then $M_{0}\cong (2^{t-1}\times \Omega_{3}(q)^t\cdot 2^{t-1})\cdot S_{t}.$ By the condition $|X|<|M_{0}|^3$, we have $q^{m^2} \prod_{j=1}^{m}(q^{2j}-1)<2^{3t-5}\cdot (t!)^3\cdot q^{3t}\cdot (q^2-1)^{3t}.$ Applying $q^{m^2}\le \prod_{j=1}^{m}(q^{2j}-1)$, we obtain that $q^{2m^2-9t}<2^{3t-5}\cdot (t!)^3.$ Note that $2m+1=3t,$ and then 
\begin{equation}\label{e3.6.2}
	q^{9t^2-24t+1}<2^{6t-10}\cdot (t!)^6.\tag{3.6.2}
\end{equation}
Assume that $t=3,$ by (\ref{e3.6.2}) we have $q^{10}<2^{14}\cdot 3^6.$ It follows that $q=3$ or $5$. However, it is easily known that $v$ is even in each case, a contradiction. If $t\ge 5,$ then $t !<5^{(t^2-3 t+1) / 3}$ by Lemma \ref{l7}(2). From (\ref{e3.6.2}), we get $q^{9t^2-24t+1}<2^{6t-10}\cdot5^{2t^2-6 t+2}<q^{3t^2-3t-7},$ and so $6t^2<21t-8$, which is impossible for $t\ge5.$ Therefore, $i\ge 5.$\\
\indent Suppose that $i\ge 5.$ Similarly, we have $q^{2m^2}<2^{6t-5}\cdot (t!)^3\cdot q^{3ti(i-1)/2},$ and so 
\begin{equation}\label{e3.6.3}
	q^{4m^2-3ti(i-1)}<2^{12t-10}\cdot (t!)^6\tag{3.6.3}
\end{equation}
If $t=3$, we have $q^{3i+1}<2^{32}\cdot 3^{6}$ from the inequality (\ref{e3.6.3}). It follows that $i=5,7.$ If $i=5,$ then $m=7$, and so $q^{16}<2^{32}\cdot 3^{6}$. The inequality holds only for $q\in\{3,5\}$. In the same way, we obtain that $q=3$ if $i=7.$ Hence, we conclude that the possible values of $(t,i,q)$ are $(3,5,3),(3,5,5)$, $(3,7,3).$ On the other hand, we observe that $v$ is even in each case, a contradiction. Thus $t\ge 5,$ in which case $q^{4m^2-3ti(i-1)}<2^{12t-10}\cdot 5^{2(t^2-3t+1)}<q^{3t^2+3t-7}$, and then $4m^2-3ti(i-1)<3t^2+3t-7.$ Note that $2m+1=it$, we have $i^2(t^2-3t)+it<3t^2+3t-8.$ It follows that $25(t^2-3t)+5t<3t^2+3t-8$ as $i\ge 5,$ i.e., $22t^2<73t-8,$ which is impossible.\vspace{2mm}\\
{\bfseries Case (3)}. Let $X=P\Omega_{7}(q)$ with $q$ prime and $q\equiv \pm3\pmod{8}$, and $M_{0}=\Omega_{7}(2)$. \\
\indent In this case, we get
\begin{equation*}
v=\dfrac{q^9(q^6-1)(q^4-1)(q^2-1)}{2^{10}\cdot 3^{4}\cdot 5\cdot 7}
\end{equation*}
and $n\le 7.$ By Lemma \ref{l1}(2), we have $q^9(q^6-1)(q^4-1)(q^2-1)<2^{11}\cdot 3^{4}\cdot 5\cdot 7^3$, which is impossible. $\hfill \square$ \vspace{3mm}
\begin{proposition}\label{prop7}
	If $\mathcal{D}$ and $G$ satisfy Hypothesis $1$, then $X$ cannot be $P\Omega_{2 m}^{\epsilon} (q),$ where $m\ge 4$ and $\epsilon=\{+,-\}$.
\end{proposition}
{\bfseries Proof}.  Let $M_0=M \cap X$, where $M=G_\alpha$ with $\alpha$ is a point of $\mathcal{P}$. By Lemma \ref{l5}, one of the following holds:
\begin{enumerate}
		\item [\rm (1)]$X=P\Omega_{8}^+(q)$, $M_{0}=\Omega_{8}^+(2)$ with $q$ prime and $q\equiv \pm3\pmod{8}$;
		\item [\rm (2)]$X=P\Omega_{8}^+(q)$, $M_{0}=2^3\cdot 2^6\cdot PSL_{3}(2)$, where $q=p\equiv \pm3\pmod{8}$, and $G$ contains a triality automorphism of $X$;
		\item [\rm (3)]$M_{0}$ is the parabolic subgroup of $X$ with $q$ even;
		\item [\rm (4)]$M$ is the stabilizer of a non-singular subspace with $q$ odd;
		\item [\rm (5)]$M_0$ is the stabilizer of an orthogonal decomposition $V=V_1 \oplus \cdots \oplus V_t$ with $\operatorname{dim}\left(V_j\right)$ is constant $i$ with $q$ odd;
		\item [\rm (6)]$q=q_0^t$ is odd with $t$ odd prime, and $M=N_G(X)$.
\end{enumerate}
Next, we will consider each cases in turn.\vspace{2mm}\\
{\bfseries Cases (1)-(2)}. Suppose that $(X,M\cap X)=(P\Omega_{8}^+(q),\Omega_{8}^+(2))$ or $(P\Omega_{8}^+(q),2^3\cdot 2^6\cdot PSL_{3}(2))$, where $q=p\equiv \pm3\pmod{8}$. If $(X,M\cap X)=(P\Omega_{8}^+(q),\Omega_{8}^+(2))$, then
\begin{center}
	$v=|X:M\cap X|=\dfrac{q^{12}(q^6-1)(q^4-1)^2(q^2-1)}{2^{14}\cdot 3^{5}\cdot 5^2\cdot 7},$ 
 \end{center}
and if $(X,M\cap X)=((P\Omega_{8}^+(q),2^3\cdot 2^6\cdot PSL_{3}(2))$, we have 
\begin{equation*}
v=\dfrac{q^{12}(q^6-1)(q^4-1)^2(q^2-1)}{2^{14}\cdot 3\cdot 7},
\end{equation*}
respectively. Moreover, we conclude that $n\le 7$ by Proposition \ref{prop2}. Then $v>2n^2\ge 2k$, a contradiction.\vspace{2mm}\\ 
{\bfseries Case (3)}. Suppose that $q$ is even, and $M_{0}$ is the parabolic subgroup of $X$. We will deal with the case where $(m,\epsilon)=(4,+)$ and $G$ contains a triality automorphism at the end of this part.\\
\indent Firstly, assume that $M$ stabilizes a totally singular $i$-space, and $i\le m-1$. If $i=m-1$ and $\epsilon=+$, then $M=P_{m,m-1}$, otherwise $M=P_{i}$. Note that 
$v>2^{-5}q^{2mi-\frac{3}{2}i^2-\frac{1}{2}i}$ from \cite[Proposition 4.4]{lprime}.
   By Proposition \ref{prop2}, we know that $n\le q^{m-1}+1.$ On the other hand, there is a subdegree $d$ of $X$ which is a power of $p$ except for the case where $\epsilon=+,$ $m$ is odd and $M=P_{m}$ or $P_{m-1}$ from Lemma \ref{l4}. Moreover, $(v-1)_{p}=q^2$ or $8$, and so $k^*\le q^3.$ Now we have
\begin{equation}\label{e3.7.1}
	q^{(4mi-3i^2-i)/2}<2^6\cdot q^3\cdot (q^{m-1}+1)<2^7\cdot q^{m+2}\tag{3.7.1}
\end{equation}
from Lemma \ref{l1}(2). By (\ref{e3.7.1}), we obtain that $4mi-3i^2-i<2m+18$, and so $2m(2i-1)-3i^2-i<18$. Note that $m\ge i+1$, then $i^2+i<20.$ Hence $i=1,2,3.$\\
\indent If $i=1$, in this case, $v=\frac{(q^m-\epsilon)(q^{m-1}+\epsilon)}{q-1}$ and $k^*\le q$. Assume that $\epsilon=+.$ Then we have $v=\frac{(q^m-1)(q^{m-1}+1)}{q-1}$, and so 
\begin{equation*}
	\frac{(q^m-1)(q^{m-1}+1)}{q-1}=v<2q(q^{m-1}+1)
\end{equation*}
from Lemma \ref{l1}(2). The previous inequality implies that $q^m-1<2q(q-1)$, contrary to $m\ge 4.$ If $\epsilon=-,$ we have $v=\frac{(q^m+1)(q^{m-1}-1)}{q-1}$, and so
$q^m\cdot q^{m-2}<\frac{(q^m+1)(q^{m-1}-1)}{q-1}<2q(q^{m-1}+1)<3q^m.$ It follows that $q^{m-2}<3,$ a contradiction.\\
\indent If $i=2$, we know that  $v=(q^{2m-2}-1)(q^m-\epsilon)(q^{m-2}+\epsilon)/(q-1)(q^2-1)$ and $k^*\le q$. Lemma \ref{l1}(2) implies that
\begin{equation*}
	\frac{(q^{2m-2}-1)(q^m-\epsilon)(q^{m-2}+\epsilon)}{(q-1)(q^2-1)}=v<2q(q^{m-1}+1),
\end{equation*}
 that is, $(q^{m-1}-1)(q^m-\epsilon)(q^{m-2}+\epsilon)<2q(q-1)(q^2-1)$. Note that $m\ge 4$, then $(q^2+q+1)(q^4+1)<2q$, which is impossible. In the same way, we assert that the remaining case $i=3$ is also impossible. \\
\indent Next, we consider that $M=P_{m}$ when $X=P\Omega_{2m}^+(q).$ Recall that here $P_{m-1}$ and $P_{m}$ are the stabilizers of totally singular $m$-spaces from the two different $X$-orbits. As $v=|X:M_{0}|=(q^{m-1}+1)(q^{m-2}+1)\cdots (q+1),$ then $v>q^\frac{m(m-1)}{2}$. Since $n$ is an odd prime divisor of $|M_{0}|_{p'},$ we conclude that $n\le \frac{q^m-1}{q-1}$. Assume that $m$ is even, then there is a subdegree $d$ which is a power of $p$. On the other hand, $q$ is the highest power of $p$ dividing $v-1$ by \cite{linear}. From Lemma \ref{l1}(4), we know that $k^*\le q.$ Since $v<2k$, we have $q^\frac{m(m-1)}{2}<2q\cdot \frac{q^m-1}{q-1}< q^{m+2},$ and then $m(m-1)<2m+4,$ which is impossible for $m\ge 4.$  If $m$ is odd, then $q^\frac{m(m-1)}{2}<\frac{2(q^m-1)^2}{(q-1)^2}<q^{2m+1}$, which implies that $m(m-1)<4m+2,$ and hence $m=5.$ Then the action here is of rank $3$, with non-trivial subdegrees $\frac{q(q^2+1)(q^5-1)}{q-1}$ and $\frac{q^6(q^5-1)}{q-1}.$ However,  Dempwolff \cite{De} has given the classification of the symmetric designs with a primitive rank 3 permutation group, none of the parameters satisfy the order $n$ is prime.\\
\indent Finally, we suppose that $(m,\epsilon)=(4,+)$ and $G$ contains a triality automorphism. 
We use the result of \cite{maximal}, where the maximal subgroups are determined. If $M\cap X$ is a parabolic subgroup of $X$, then it is either $P_{2}$ or $P_{134}$. The former has been rule out above, so consider the latter. If $M\cap X=P_{134},$ then $v=|X:M\cap X|=\frac{(q^6-1)(q^4-1)}{(q-1)^3}>q^{11}.$
 By Proposition \ref{prop2}, $n$ is an odd prime divisor of $|M_{0}|_{p'}$, that is, $n\le q+1.$ Hence $q^{11}<2(q+1)^2$, a contradiction.\vspace{2mm}\\
{\bfseries Case (4)}. Suppose that $q$ is odd, and $M$ is the stabilizer of a non-singular subspace. Here $M_{0}=N_{i}$ with $i\le m$ and $s=2m.$\\
\indent If $i=1$, then $M_{0}\cong \hat{ }\Omega_{2m-1}(q)\cdot 4$ , and so $v=|X:M\cap X|=\frac{q^{m-1}(q^m-\epsilon)}{2}.$ By Proposition \ref{prop2}, then $n\le \frac{q^{m-1}+1}{2}.$ Moreover, we get $k^*\le q$ from \cite{R3}. Then Lemma \ref{l1}(2) implies that $\frac{q^{m-1}(q^m-\epsilon)}{2}=v<q(q^{m-1}+1)$. Note that $q^m-\epsilon>q^{m-1}+1,$ then $q^{m-1}<2q$, and so $q^2\le q^{m-2}<2$, which is impossible. Therefore, $2\le i\le m.$\\
\indent Assume that $2\le i\le m.$ Then by \cite[Proposition 4.4]{lprime}, we know $v>\frac{q^{i(s-i)}}{4}.$ By Proposition \ref{prop2}, we know that $n$ is an odd prime divisor of $|M_{0}|_{p'},$ and so $n\le \frac{q^{(s-i)/2}+1}{2}$. Since $v<2k,$ then 
\begin{equation*}
	\dfrac{q^{i(s-i)}}{4}<v<2n^2\le \dfrac{(q^\frac{s-i}{2}+1)^2}{2},
\end{equation*}
which implies that $q^{i(s-i)}<q^{s-i+2}$, that is, $s(i-1)<i^2-i+2.$ Recall that $s=2m\ge 2i,$ and then $i^2<i+2,$ a contradiction.\vspace{2mm}\\
{\bfseries Case (5)}. Assume that $q$ is odd, and $M_{0}$ is the stabilizer of an orthogonal decomposition $V=V_1 \oplus \cdots \oplus V_t$ with $\operatorname{dim}\left(V_j\right)=i$, where $2m=it,t\ge 2$.\\
\indent Let $i=1.$ By Proposition \ref{prop2}, we have $|X|<2|M\cap X|\cdot |M\cap X|_{p'}^2$, and so the pairs satisfy this condition are $(X,M_{0})=(P\Omega_{8}^+(3),2^6\cdot A_{8})$ and $(P\Omega_{10}^-(3),2^8\cdot A_{10})$. Note that $n\le 7$ in these cases. For all of these possibilities, $v>2n^2\ge 2k$, a contradiction. Hence $i\ge 2.$\\
\indent Assume first that $i$ is even, then $M_{0}\cong d^{-1}\Omega_{i}^{\epsilon_{1}}(q)^t\cdot 2^{2(t-1)}.S_{t},$ where $\epsilon=\epsilon_{1}^t$ and $d\in \{1,2,4\}.$ If $i=2$, then $v>\frac{q^{2t^2-t}}{2^{t-2}\cdot t!\cdot (q+1)^t}$ from Lemma \ref{l8}. By Proposition \ref{prop2}, $n\le \frac{t(q+1)}{2},$ and so by Lemma \ref{l1}(2), we conclude that
\begin{equation*}
	\dfrac{q^{2t^2-t}}{2^{t-2}\cdot t!\cdot (q+1)^t}<v<2n^2\le  \dfrac{t^2(q+1)^2}{2}.
\end{equation*}
 Simple calculation gives $q^{2t^2-2t-2}<2^{2t-1}\cdot t!\cdot t^2.$ Note that $t=m$ and $m\ge 4$, then  $t !<2^{4 t(t-3)/3}$ by Lemma \ref{l7}. On the other hand, we know that $t^6\le 3^{2t},$ and hence $q^{6t^2-6t-6}<2^{6t-3}\cdot2^{4 t(t-3)}\cdot 3^{2t}<q^{4t^2-4t-3}$. It follows that $2t^2<2t+3,$ which is impossible.\\
\indent If $i=4,$ then $m=2t$. From Proposition \ref{prop2}, we know that $|X|<|M_{0}|^3,$ and so
\begin{equation}\label{e3.7.2}
	q^{m(2m-1)}<2^{6t-3}\cdot (t!)^3\cdot q^{6t} (q^4-1)^{3t}<2^{6t-3}\cdot (t!)^3\cdot q^{18t}.\tag{3.7.2}
\end{equation}
If $t=2$, then $m=2t=4$, and $v=\frac{q^8(q^2+\epsilon_{1})(q^6-1)(q^2+1)}{2}$ is even, a contradiction. If $t=3$, then $m=6$ and $q^{12}<2^{18}\cdot 3^3$ from (\ref{e3.7.2}). The previous inequality yields $q=3.$ So $(X,M_{0})=(P\Omega_{12}^+(3), \hat{ }\Omega_{4}^+(3)^3.2^2.S_{3})$ and $v>2\cdot 3^2\ge 2\cdot n^2\ge 2k$, a contradiction. If $t\ge 4$, then $q^{8t^2-20t}<2^{4t^2-6t-3}$ by Lemma \ref{l7}. It implies that $4t^2<14t-3,$ which is impossible.\\
\indent If $i\ge 6$, by Proposition \ref{prop2}, we have $q^{m(2m-1)}<2^{6t-3}\cdot (t!)^3\cdot q^\frac{3it(i-1)}{2}$. Note that $2m=it$, then 
\begin{equation}\label{e3.7.3}
	q^{i^2t^2+2it-3i^2t}<2^{12t-6}\cdot (t!)^6.\tag{3.7.3}
\end{equation}  
Similarly, if $t=2$, let $i=2a$ for some integer $a$, then we get $v$ is even, a contradiction. If $t=3$, then $q^{6i}<2^{36}\cdot 3^6$ by (\ref{e3.7.3}). Since $q\ge 3$, we have $3^{6i}<2^{36}\cdot 3^6$, which conflicts that $i\ge 6.$  If $t\ge 4$, then $q^{i^2t^2+2it-3i^2t}<2^{8t^2-12t-6}$ by Lemma \ref{l7}. Then $(i^2t^2+2it-3i^2t)\cdot \text{log}_{p}q<(8t^2-12t-6)\cdot \text{log}_{p}2<(8t^2-12t-6)\times 0.7$, and so $10(i^2t^2+2it-3i^2t)<56t^2-84t-42$. Note that $i\ge 6$, then $152t^2<438t-21$, contrary to $t\ge 4.$\\
\indent Next, we assume that $i$ is odd. Then by \cite[Proposition 4.2.14]{classical}, we know that $M_{0}\cong (2^{t-2}\times \Omega_{i}(q)^t.2^{t-1}).S_{t}$ with $t$ even and $\epsilon=(-1)^\frac{m(q-1)}{2}.$ If $i=3$, then $3t=2m$ and $t\ge 4.$ By Proposition \ref{prop2}, we have 
\begin{equation*}
	q^{m(2m-1)}<2^{6t-6}\cdot (t!)^3\cdot q^{3t} (q^2-1)^{3t}\cdot 2^{-3t},
\end{equation*}
which yields that $q^{9t^2-3t}<2^{6t-12}\cdot (t!)^6\cdot q^{18t}.$ Note that $t\ge 4$, then $q^{9t^2-21t}<2^{8t^2-18t-12}$ by Lemma \ref{l7}. It follows that $9t^2-21t<8t^2-18t-12,$ and so $t^2<3t-12$, a contradiction. Thus $i\ge 5.$ Here if $t=2$, then $m=i$ and $i=2a+1$ for some integer $a$. In this case,
\begin{equation*}
	v=\dfrac{q^{3a^2+2a}(q^{2a+1}-\epsilon_{1})(q^{4a}-1)\cdots (q^2-1)}{2(q^{2a}-1)^2\cdots (q^2-1)^2}
\end{equation*}
 is even, a contradiction. Therefore, $t\ge 4.$ By Proposition \ref{prop2}, $|X|<|M_{0}|^3$, we get  $q^{m(2m-1)}<2^{6t-6}\cdot (t!)^3\cdot q^\frac{3it(i-1)}{2}$. Lemma \ref{l7}(1) implies that $q^{i^2t^2+2it-3i^2t}<2^{8t^2-12t-12},$ and so $i^2(t^2-3t)+i\cdot 2t<8t^2-12t-12,$ which is impossible for $i\ge 5$.\vspace{2mm}\\
{\bfseries Case (6)}. Finally, we assume that $q=q_0^t$ is odd with $t$ odd prime, and $M=N_G(X)$.\\
\indent By \cite[Proposition 4.5.10]{classical}, the subgroup $M_{0}\cong P\Omega_{2m}^\epsilon(q_{0})$ with $m\ge 4.$ If $t=3$, Lemma \ref{l8} implies that 
$v>q_{0}^{4m^2-2m-4}.$ By Proposition \ref{prop2}, we have $n\le q_{0}^m+1.$ Moreover,
	$q_{0}^{4m^2-2m-4}<v<2n^2\le 2(q_{0}^m+1)^2<2^3\cdot q_{0}^{2m},$
and so $4m^2-2m-4<2m+3.$ It follows that $4m^2<4m+7,$ a contradiction. Thus $t\ge 5,$ combining Lemma \ref{l8} with Proposition \ref{prop2}, we have
\begin{equation*}
	q_{0}^{tm(2m-1)}<2^3\cdot q_{0}^{3m(2m-1)}\cdot (1+q_{0}^{-m})^3<2^6\cdot q_{0}^{6m^2-3m}. 
\end{equation*}
 Note that $t\ge 5,$ then $q_{0}^{4m^2-2m}<2^6,$ and so $2m^2<m+3$. The desired contradiction.\vspace{3mm}$\hfill \square$
\subsection{Alternating groups }  
\indent Let $\mathcal{D}$ be a 2-$(v,k,\lambda)$ non-trivial symmetric design of prime order, and $G\le Aut(\mathcal{D})$ be flag-transitive with $X=Soc(G)=A_{m}$ $(m\ge 5)$. Recall that if $m=6$, then $G\cong A_{6}$$,S_{6},M_{10},PGL(2,9)$ or $P\Gamma L(2,9)$, and if $m\neq 6$, then $G\cong A_{m}$ or $S_{m}$.\\
\indent We first consider the case $m=6$ and $G\cong M_{10},PGL(2,9)$ or $P\Gamma L(2,9)$. For each group $ M_{10},PGL(2,9)$ or $P\Gamma L(2,9)$, we know their indices of maximal subgroups are $2,10,36$ and $45$, so $v=10,36$ or $45$. If $v=10$, we conclude that $k^*=3$ or $9$ by Lemma \ref{l1}(4), and $k=3n$ or $9n$, respectively. Note that $k\mid |G_{\alpha}|$, if $G=M_{10}$ or $PGL(2,9)$, then $k\mid 72$ and if $G=P\Gamma L(2,9)$, then $k\mid 144$. Now we get the possible parameters $(v,k)=(10,6)$ or $(10,9)$. However, the parameters $(v,k,\lambda)=(10,9,8)$ is trivial. Then $k=6$. From Lemma \ref{l1}(2), we have $\lambda=\frac{k(k-1)}{v-1}=\frac{10}{3},$ a contradiction. If $v=36$, we have $k^*=5,7$ or $35$ by the condition $k^*\mid v-1$. Since $k\mid |G_{\alpha}|$, then $k^*=5$ and $n=2$. But $\lambda=\frac{k(k-1)}{v-1}=\frac{18}{7}$ is not integer, a contradiction. Similarly, if $v=45$, from Lemma \ref{l1}(4) and $k^*\mid |G_{\alpha}|$, we obtain that $k^*=4$. Further, we have $\lambda=\frac{14}{11}$, a contradiction. Therefore, if $m=6$, then $G\cong A_{6},S_{6}$.\\
\indent Let $G= A_{m}$ or $S_{m}$ with $m\ge 5$ in the following. Assume that $\alpha \in \mathcal{P}$, then $G_{\alpha}$ acts on $\mathcal{P}$, and it also acts on the set $\Omega=\{1,2,\cdots,m\}$. According to Lemma \ref{l6}, we consider the following three cases: (1) $M:=G_{\alpha}$ is intransitive on $\Omega$;
(2) $M:=G_{\alpha}$ is transitive and imprimitive on $\Omega$;
(3) $M:=G_{\alpha}$ is primitive on $\Omega$.  
\begin{proposition}\label{prop8}
	If $\mathcal{D}$ and $G$ satisfy Hypothesis $1$, then $X\neq A_{m}$ $(m\ge 5).$
\end{proposition}
{\bfseries Proof}. Let $G=S_{m}$ or $A_{m}$ where $m\ge 5.$\vspace{2mm}\\
{\bfseries Case (1)}. Suppose that $M:=G_{\alpha}$ is intransitive on $\Omega$, then $M=(S_{s}\times S_{t})\cap G$, where $m=s+t$ and  $s\neq t$. Without loss of generality, assume that $s<\frac{m}{2}$. Since $G$ is flag-transitive, then $M$ is transitive on the blocks through $\alpha$. Furthermore, $M$ fixes exactly one point $\alpha$, and stabilizes only one $s$-subset $\Theta\subseteq \Omega$,
hence we have the points of $\mathcal{P} $ corresponds to all the set of $s$-subset of $\Omega$. So that $v=\binom{m}{s}$ and $G$ acting on $\mathcal{P}$ has rank $s+1$. The subdegrees of $G$ are:
\begin{equation}\label{e3.2.1}
	d_{0}=1, d_{i+1}=\binom{s}{i}\binom{m-s}{s-i},i=0,\dots ,s-1.\tag{3.2.1}
\end{equation}
\indent Since $k^* \mid d$ for any subdegree $d$ of $G$ and $d_{s}=s(m-s)$ is a subdegree of $G$, then $k^* \mid s(m-s)$ and $k^*\le s(m-s)<(m-s)^2$. From Lemma \ref{l1}(3), we know that $n$ is a prime factor of $|G_{\alpha}|,$ and so $n\le (m-s)$. Lemma \ref{l1}(2) implies that $\binom{m}{s}=v<2k=2nk^*\le 2(m-s)^3$.\\
\indent We first assert that $s\le 4$. In fact, if $s\ge 6$, then $6\le s<t$, and hence $t\ge 7$. Now
\begin{equation*}
	\dfrac{\binom{m}{s}}{2(m-s)^3}>\dfrac{2^{s-5}\binom{m-(s-5)}{5}}{2(m-s)^3}=\dfrac{2^{s-6}\binom{t+5}{5}}{t^3}\ge \dfrac{\binom{t+5}{5}}{t^3}=\dfrac{(t+5)\cdots (t+1)}{120t^3}>1,
\end{equation*}
which implies that $\binom{m}{s}>2(m-s)^3$, contrary to $v<2k$. So $s\le 5$. If $s=5$, then $\binom{m}{5}< 2(m-5)^3$, which is impossible for $m\ge 11$. Thus, $s=1,2,3$ or $4$.\\
\indent If $s=1$, then $v=m$ and $G$ acts on $\mathcal{P}$ as an alternating or symmetric group by its natural action. The group $G$ is $(m-2)$-transitive. Note that $2<k<v-1=m-1$, then $G$ is $k$-transitive on $\mathcal{P}$. We have $|\mathcal{B}|=|B^G|=\binom{m}{k}=m$ for any block of $\mathcal{B}$, and so $k=1$ or $m-1$. Hence $\mathcal{D}$ is trivial.\\
\indent If $s=2$, from (\ref{e3.2.1}), we get the subdegrees are $d_{0}=1,d_{1}=\binom{m-2}{2}$ and $d_{2}=2(m-2)$. Then $G$ is a primitive rank 3 group acting on the points of $\mathcal{D}$. By \cite[Theorem]{De}, we obtain that $(v,k,\lambda)=(15,7,3)$ with $Soc(G)=A_{6}$, or $(v,k,\lambda)=(35,17,8)$ with $Soc(G)=A_{8}.$
But none of these parameters satisfy the condition $n=k-\lambda$ is prime.\\
\indent If $s=3,$ by the inequation $\binom{m}{3}<6(m-3)^2$, we have $7\le m \le 32$. By the same way, if $s=4,$ we have $9\le m \le 12$. Now we compute all possible parameters $(v,k,\lambda)$ satisfying the following conditions:
\begin{enumerate}
	\item[\rm (a)]$v=\binom{m}{s}$;
	\item[\rm (b)]$k^*\mid (v-1,d_{s})$, where $d_{s}=s(m-s)$ is a subdegree of $G$;
	\item[\rm (c)]$n\le (m-s)$ and $n$ is a prime;
	\item[\rm (d)]$vk\mid m!$;
	\item[\rm (e)]$\lambda=\dfrac{k(k-1)}{v-1}.$
\end{enumerate}
With the aid of the computer algebra system \cite{gap}, we obtain that the possible 5-tuples $(v,k,\lambda,s,m)=(120,35,10,3,10),(1771,60,2,3,23),(4495,322,23,3,31)$. However, $k-\lambda$ is not prime for each possible care.\vspace{2mm}\\
{\bfseries Case (2)}. Assume that $G_{\alpha}$ is transitive and imprimitive on $\Omega$, then $G_{\alpha}=(S_{s}\wr S_{t})\cap G$, where $m=st$. Let $\Sigma=\{\Delta_{0},\Delta_{1},\dots ,\Delta_{t-1}\}$ be a non-trivial partition of $\Omega$ preserved by $G_{\alpha}$, where $|\Delta_{i}|=s$, $i \in\{0,1,\dots,t-1\}$ and $s,t\ge 2$. Since $G$ is flag-transitive, then $G_{\alpha}=G_{\Sigma}$. By the argument in \cite[Section 3]{finite}, we identify $\mathcal{P}$ with all $t\times s$ imprimitive partitions of $\Omega$. Hence, 
\begin{equation}\label{e3.2.2}
	v=\dfrac{(ts) !}{(s !)^t t !}=\dfrac{\binom{ts}{s}\cdots\binom{3s}{s}\binom{2s}{s}}{t!},\tag{3.2.2}
\end{equation}
and the subdegrees of $G$ are
\begin{equation*}
	\begin{cases}
		d_{j}=\frac{1}{2}\binom{t}{j}\binom{s}{1}^j=2^{j-1}\binom{t}{j} ,  &      s=2,\\
		d_{j}=\binom{t}{j}\binom{s}{1}^j=s^j\binom{t}{j},    &s\ge 3.
	\end{cases}
\end{equation*}
\indent If $s=2$, then $v=(2t-1)(2t-3)\cdots 5\cdot 3$ and $d_{2}=t(t-1)$. From Lemma \ref{l1}, we have $k=nk^*\le t^2(t-1)$. Recall that $v<2k$, then 
\begin{equation*}
	(2t-1)(2t-3)\cdots 5\cdot 3<2t^2(t-1),
\end{equation*}
which implies that $t\le 3$. Since $m=st=2t$ and $m\ge 5$, then $t\ge 3$. So $t=3$. By straightforward calculation, we obtain that $v=15$ and $d_{2}=6$. Combining it with $n\le t$ and $k^*\mid d_{2}$, then $k=4,6,9$ or $12$. But the value of $\lambda$ is not integer, a contradiction. Thus, $s\ge 3$.\\
\indent If $s\ge 3$, then $d_{2}=s^2\binom{t}{2}$. Since 
\begin{equation*}
	\dfrac{\binom{js}{s}}{j}=\binom{js-1}{s-1}=\dfrac{js-1}{s-1}\cdot \dfrac{js-2}{s-2}\cdots \frac{js-(s-1)}{1}>j^{s-1},
\end{equation*}
we get $v>(t!)^{s-1}$ from (\ref{e3.2.2}). Recall that $n$ is a prime divisor of $|M|$ and $k^*\mid d_{2}$, then $k=k^*n\le s^2\binom{t}{2}\cdot st= s^3t\binom{t}{2}$. Now we have 
\begin{equation*}
	(t!)^{s-1}<v<2k\le 2s^3t\binom{t}{2}=s^3t^2(t-1).
\end{equation*}
It follows that (i) $t=2$ and $3\le s\le 14$, (ii) $t=3$ and $3\le s\le 5$, or (iii) $t=4$ and $s=3$.
For each of the possible 2-tuples $(t,s)$, we have the value of 2-tuples $(v,d_{2})$. Note that $n\le \text{max}\{s,t\}$, then $v<2\cdot \text{max}\{s,t\}\cdot d_{2}$. By simple calculation, we know the possible 4-tuples $(t,s,v,m)=(2,3,10,6),(2,4,35,8)$ or $(2,5,126,10)$. Furthermore, the parameters $(v,k,\lambda)$ satisfy the following arithmetical conditions:
\begin{enumerate}
	\item [\rm (a)]$vk\mid |G|$, where $G=A_{m}$ or $S_{m}$;
	\item [\rm (b)]the basic equation $\lambda=\dfrac{k(k-1)}{v-1}$;
	\item [\rm (c)]the order $n=k-\lambda$ is prime.
\end{enumerate}
By using GAP \cite{gap}, we found that no parameter satiafies the conditions (a)-(c). Thus, $G_{\alpha}$ cannot act transitively and imprimitively on $\Omega$.\vspace{2mm}\\
{\bfseries Case (3)}. From Lemma \ref{l6}, we know $v=15$ and $G=A_{7}$ or $A_{8}$. Note that $k^*\mid v-1,$ then $k^*=2,7$ or $14.$ Assume that $k^*=2$, then $k=2n$. Lemma \ref{l1}(3) implies that $2n\mid 168$ or $1344$. Note that $n$ is a prime number, then $n=2,3$ or $7.$ Next, for $k^*=7$ or $14$, it is easy to know that $(v,k)=(15,4),(15,6) $ or $(15,14).$ The first two cases are excluded since $\lambda=\frac{k(k-1)}{v-1}$ is not an integer, and the last case is trivial. $\hfill \square$ \vspace{3mm}\\
{\bfseries Proof of Theorem \ref{th1}}.  
By Propositions \ref{prop1} and \ref{prop8}, we conclude that $X$ cannot be a sporadic simple group, a finite simple exceptional group and an alternating group. Assume that $\mathcal{D}$ is a
symmetric design of prime order, and  admitting flag-transitive automorphism group $G$ with $X$ is a classical group.   Since the order $n=k-\lambda$ is prime, then $(k,\lambda)=1$ or $n$. If $(k,\lambda)=1$, by \cite{1}, we know that the parameters of $\mathcal{D}$ is $(v,k,\lambda)=(11,5,2)$ with $G=PSL_{2}(11)$ or $(v,k,\lambda)=(n^2+n+1,n+1,1)$ with $X=PSL_{3}(n)$. If $(k,\lambda)=n,$ we may assume that $k^*\ge 3$ by the main theorem in \cite{lprime}. Moreover, 
If $X=PSL_{2}(q)$, then $\mathcal{D}$ either is a $2$-$(7,4,2)$ design or $2$-$(11,6,3)$ design in \cite{psl}. By \cite{four}, we know that if $X=PSL_{m}(q)$ and $m\le 4$, except for the above possibilities, then $\mathcal{D}$ is  $2$-$(n^2+n+1,n^2,n^2-n)$ design, that is, the complement of $PG_{2}(n),$ and $PSL_{3}(n)\unlhd G\le P\Gamma L_{3}(n)$. Hence, we focus on the dimension $m\ge 5$ when deal with $X$ is a linear or unitary group. Now Theorem \ref{th1} follows from Proposition \ref{prop3}-\ref{prop7}.$\hfill \square$ \vspace{3mm}\\

\end{document}